\documentstyle[12pt,fleqn]{article}
\newcommand{\tn}{\otimes}

\newcommand{\g}{\mbox{\bf g}}
\renewcommand{\t}{\mbox{\boldmath $t$}}
\newcommand{\tx}{\mbox{\footnotesize\boldmath $t$}}

\newcommand{\id}{\mbox{\rm id}}
\newcommand{\qq}{{\cal Q}}
\newcommand{\rr}{{\cal R}}
\newcommand{\ff}{{\cal F}}

\newcommand{\eins}{1\hspace{-3pt}{\rm I}}
\newcommand{\cc}{{\rm C}\!\!\!\mbox{\small\sf l}\;}

\newcommand{\beq}{\begin{equation}}
\newcommand{\eeq}{\end{equation}}
\newcommand{\bgar}{\begin{array}}
\newcommand{\ear}{\end{array}}
\newcommand{\con}{\begin{array}[b]{c}\mbox{\footnotesize $\;\leftharpoonup$}\\}
\newcommand{\ccon}{\begin{array}[b]{c}\mbox{\footnotesize $\;\rightharpoonup$}
\\}
\newcommand{\crr}{\mbox{\asn
$\!\!\begin{array}[b]{c}\mbox{\footnotesize $\:\!\rightharpoonup$}\\ \rr\ear
\!\!$\ase}}
\newcommand{\ccf}{\mbox{\asn
$\!\!\begin{array}[b]{c}\mbox{\footnotesize $\,\,\leftharpoonup$}\\ \ff\ear
\!\!$\ase}}
\newcommand{\cf}{\mbox{\asn
$\!\!\begin{array}[b]{c}\mbox{\footnotesize $\,\,\rightharpoonup$}\\ \ff\ear
\!\!$\ase}}
\newcommand{\ccv}{\mbox{\asn
$\!\!\begin{array}[b]{c}\mbox{\footnotesize $\:\!\leftharpoonup$}\\ v\ear\!\!$
\ase}}
\newcommand{\cv}{\mbox{\asn
$\!\!\begin{array}[b]{c}\mbox{\footnotesize $\:\!\rightharpoonup$}\\ v\ear\!\!$
\ase}}
\newcommand{\ccz}{\mbox{\asn
$\!\!\begin{array}[b]{c}\mbox{\footnotesize $\;\!\leftharpoonup$}\\ z\ear\!\!$
\ase}}
\newcommand{\cz}{\mbox{\asn
$\!\!\begin{array}[b]{c}\mbox{\footnotesize $\,\rightharpoonup$}\\ z\ear\!\!$
\ase}}

\begin{document}

\newcommand{\ase}{\renewcommand{\arraystretch}{1.4}}
\newcommand{\asn}{\renewcommand{\arraystretch}{0}}

\ase
\thispagestyle{empty}
\begin{flushright}
LMU--TPW 95--13 \\ August 1995 \\ q-alg/9509001
\end{flushright}
\vskip 1cm
\bigskip\bigskip
\begin{center}
{\bf\LARGE{On the Drinfel'd-Kohno Equivalence}}
\end{center}

\begin{center}
{\bf\LARGE{of Groups and Quantum Groups}}
\end{center}

\vskip 1.5truecm
\begin{center}{\bf Ralf A. Engeldinger}\footnotetext{Supported by
Graduiertenkolleg {\sc Mathematik im Bereich ihrer Wechsel\-wirkung
\mbox{\sc mit
der Physik}} at Mathematisches Institut der Ludwig-Maximilians-Universit\"at
M\"unchen}
\\
\end{center}

\begin{center}
{\bf }
\vskip10mm
Sektion Physik der Ludwig-Maximilians-Universit\"at M\"unchen\\
Lehrstuhl Professor Wess\\ Theresienstra\ss e 37, D-80333 M\"unchen,
Federal Republic of Germany\\

\vskip5mm
e-mail: engeldin@lswes8.ls-wess.physik.uni-muenchen.de
\end{center}
\vskip 1cm
\nopagebreak
\begin{abstract}
A method to calculate matrix representations of the twist element $\ff$
of Drinfel'd -- chosen to be unitary -- is given
and illustrated at some examples. It is observed that for these F-matrices the
crystal
limit $q\!\to\! 0$ exists and that \mbox{F-matrices} twisting from $0$ to $q$
are of a
simpler form
than F-matrices twisting from $1$ to $q$. These results lead to a new
interpretation of
$q$-deformation in terms of tensor products of finite-dimensional
representations of
compact simple Lie groups.
\end{abstract}

\newpage
\section{Introduction}

Drinfel'd published his ground breaking work on deformation of semi-simple Lie
algebras
in 1989 in Russian and in 1990 in English [1]. Combining ideas from areas of
physics and
mathematics as diverse as Lie theory, Hopf algebra theory, integrable models,
conformal
field theory, gauge theory, cohomology theory and category theory he proved
their existence
and uniqueness using the mathematical frame of
{\sl quasitriangular quasi-Hopf quantized universal enveloping algebras} and
{\sl twists} as equivalence mappings between them. Interestingly he solved a
{\sl mathematical} problem by translating it into a {\sl physical} model,
thereby somehow
interchanging the roles usually played by the two sciences.

Impressive as it is the Drinfel'd-Kohno theorem has been suffering from one
crucial
drawback:
Even though the existence of the twist between the ``undeformed'' and the
``deformed'' universal enveloping algebras is proven any attempt to gain an
explicit
expression for
the twist element $\ff$ failed. It is mostly due to this circumstance that the
theorem didn't have a major impact on research papers published ever since. Its
content remained inaccesible for purposes of explicit calculations.

The present work is an attempt to overcome these difficulties.

\section{The Drinfel'd-Kohno Theorem}

{\bf Theorem }{\sc (Drinfel'd-Kohno,[1]) }
{\sl Let $\g$  be a simple finite-dimensional Lie algebra over $\cc$ and $\t$
a symmetric element of $\g\tn\g$ s.t. \mbox{$[\t,\Delta(\g)]=0$.} From this
can be constructed, on the one hand, a quasitriangular quasi-Hopf
quantized universal enveloping algebra (QtQHQUEA)
$(U\g[[h]],\Delta,\Phi,\rr)$, where $\Delta$ is the undeformed
cocommutative and coassociative coproduct and $\rr=\exp(\frac{h}{2}\t)$,
and, on the other hand, a quasitriangular Hopf quantized universal
enveloping algebra (QtHQUEA) $(U_h\g,\bar\Delta,\bar\rr)$. They are
twist-equivalent as quasitriangular quasi-Hopf algebras (QtQHAs).

{\samepage \parskip.8ex
That is there exists an invertible $\ff\in U\g^{\tn 2}[[h]]$ s.t.\footnote
{For $x=a\tn b$ let $x_{12}=a\tn b$, $x_{21}=b\tn a$ resp.
$x_{12}=a\tn b\tn\eins$,
$x_{23}=\eins\tn a\tn b$, etc.}
\beq\label{twist} \bar\Delta(x)=\ff\,\Delta(x)\,\ff^{-1} \quad
\forall x\in U\g[[h]]\eeq
\beq \bar\rr_{12}=\ff_{21}\,\rr_{12}\,\ff^{-1}_{12} \eeq
\beq\label{Phi} \eins^{\tn 3}=\bar\Phi=\ff_{23}\;(\id\tn\Delta)(\ff)\;\Phi\;
[\ff_{12}\;(\Delta\tn\id)(\ff)]^{-1}. \eeq}}

\parskip2.5ex plus0.5ex minus0.5ex

To begin with we discuss some implications:

The theorem establishes an {\sl invertible} transition from the `undeformed'
to the
`deformed' case without refering to the limit $h\!\to\! 0$. This is in clear
contrast to the
relationship between classical and quantum physics: As elements of a
commutative algebra
phase space functions have completely `forgotten' about the possibility of
being
non-commutative. To restore this information it takes an additional structure,
the Poisson
bracket, which is
independent of the commutative multiplication. To put it another way:
The commutative algebra of classical phase space functions {\sl without}
Poisson structure
contains strictly less information than the non-commutative algebra of
operators on the
Hilbert space of states.

In the case of (commutative) function algebras on groups and non-commutative
`function
algebras on quantum groups' (which are the respective dual Hopf algebras of
the universal
enveloping algebras featured in the theorem) the two different algebra
structures are
related by some kind of conjugation (or, more precisely, the dual version of
it,
`co-conjugation') or `gauge transformation' (Drinfel'd) which is, and that is
the crucial
point, invertible. This implies that both structures contain the {\sl same}
information.

Whereas we think of (non-commuting) operators and (commuting) phase space
functions as
very different kinds of objects the Drinfel'd-Kohno theorem strongly suggests
to think
of (non-commuting) `functions on a quantum group' and (commuting) functions on
a group
as the same kind of objects, and, in particular, {\sl not as primarily
linked by the limit} $h\!\to\! 0$ -- which is not to say that they are not
{\sl additionally}
linked by that limit. Their relationship is very similiar to that between a
covariant
derivative and an ordinary one. Although one might think of the latter as the
former
in the limit of vanishing coupling constant this is not a very natural way to
look at it.

\section{F-matrices}
In addition to what is stated in the Drinfel'd-Kohno theorem we will use the
following input
to explicitly calculate matrix representations of $\ff$:
\begin{enumerate}
\item The matrix representations of $\rr$ and $\bar\rr$ can be
calculated.\footnote
{The latter ones follow from the defining representations given in [2] for the
series
$A,B,C,D$.}
\item $\ff$ commutes with the coproducts of the generators of the Cartan
subalgebra of $\g$.
\item The quantized universal enveloping algebras of compact semi-simple Lie
groups are for
real $h$ equipped with a $\ast$-structure such that the matrix representations
are
$\ast$-representations.
\item $\lim_{h\to 0}\limits\ff=\eins^{\tn 2}$
\end{enumerate}

We will use the following notation: Lower case Greek letters label matrices
as a
whole, and if the same Greek label appears repeatedly in one expression we
understand this
as ordinary matrix multiplication. Thus,
 $a_\alpha b_\alpha$ means the same as $ a^i_j b^j_k$, and
 $a_\alpha b_\beta c_\alpha d_\beta$ means the same as
$a^i_j b^k_l c^j_m d^l_n$.
If a symbol has several Greek labels each of them stands for a pair of an
upper and a lower
index:
 $a_{\alpha\beta} b_{\alpha\gamma} c_\beta$ means the same as
$a^{ik}_{jl} b^{jp}_{mq}
c^l_n$.
Symbols without Greek labels will not have any matrix indices. This allows us
to denote
matrix representations of (quasi-)Hopf algebra elements by the same symbol as
the
(quasi-)Hopf algebra elements themselves with Greek labels attached to them:
$x_\alpha=\varrho_\alpha(x),\quad \rr_{\alpha\beta}
=(\varrho_\alpha\tn\varrho_\beta)(\rr)$.
Since only coinciding labels indicate index contraction, different labels may
refer
to different matrix representations (of possibly different dimensions).

\subsection{Unitarity}

{}From $\rr^\ast_{12}=\rr_{12}=\rr_{21}$ and $\bar\rr^\ast_{12}=\bar\rr_{21}$
we deduce the
identity
\beq
\ff_{12}^{-1\ast}\,\rr_{12}\,\ff_{21}^\ast
=\ff_{12}^{-1\ast}\,\rr^\ast_{12}\,\ff_{21}^\ast
=\bar\rr^\ast_{12}=\bar\rr_{21}=\ff_{12}\,\rr_{12}\,\ff^{-1}_{21}, \eeq

thus

\[ \ff^\ast\,\ff\,\rr^2=\rr^2\,\ff^\ast\,\ff. \]

{}From $\Delta(x^\ast)=\Delta(x)^\ast$ and $\bar\Delta(x^\ast)
=\bar\Delta(x)^\ast$
it follows
that
\[ \ff^\ast\,\ff\,\Delta(x)=\Delta(x)\,\ff^\ast\,\ff. \]
Thus twisting by $\ff^\ast\,\ff$ leaves the undeformed universal enveloping
algebra
invariant such that twisting by the unitary
\[ \tilde\ff=\ff (\ff^\ast\,\ff)^{-\frac12} \]
gives the same result as twisting by $\ff$. This has been noted by
Jur\v co [3].
We will
therefore assume from here without
loss of generality that $\ff$ is unitary: $\ff^\ast=\ff^{-1}$.

Since all $\rr_{\alpha\beta}$ and $\bar\rr_{\alpha\beta}$ are real matrices the
corresponding $\ff_{\alpha\beta}$ are real, orthogonal matrices:
\beq \ff_{\alpha\beta}^\top=\ff_{\alpha\beta}^{-1}. \eeq

\subsection{Orthogonal Projectors}

Denote by $H_i$, $(i=1,\ldots,$ rank $\g)$, the generators of the Cartan
subalgebra of $\g$.
In the representation $\varrho_\alpha\tn \varrho_\beta$ we write down their
expansion in terms of eigenvalues and orthonormal projectors.
\[ \Delta(H_i)_{\alpha\beta}=(H_i\tn\eins+\eins\tn H_i)_{\alpha\beta}
=\sum_c\eta_{i,c}P^{\langle i,c\rangle}_{\alpha\beta} \]
\[ P^{\langle i,c\rangle}_{\alpha\beta}P^{\langle i,c'\rangle}_{\alpha\beta}
=\delta_{cc'}P^{\langle i,c\rangle}_{\alpha\beta} \]
{}From $\Delta(H_i)\ff=\ff\Delta(H_i)$ it follows that
\beq P^{\langle i,c\rangle}_{\alpha\beta}\ff_{\alpha\beta}
=\ff_{\alpha\beta}P^{\langle i,c\rangle}_{\alpha\beta}\quad \forall i,c \eeq

Now set $\qq=\rr^2=\;$e${}^{h\tx}$ and
$\bar\qq=\bar\rr_{21}\bar\rr_{12}$,\footnote
{We use representations of $\bar\rr_{21}\bar\rr_{12}$, rather than the
conventional
$\hat R$ which only works if $\varrho_\alpha\!=\!\varrho_\beta$ (in which case
$Q\!=\!\hat R^2$).}
such that
$\bar\qq\ff=\ff\qq$.
Then, analogously, we obtain for $\qq$ and $\bar\qq$:
\[ \qq_{\alpha\beta}=\sum_k \lambda_k\,P^{(k)}_{\alpha\beta} \]
\[ \bar\qq_{\alpha\beta}=\sum_k \lambda_k\bar P^{(k)}_{\alpha\beta} \]
with
\[ P^{(k)}_{\alpha\beta}P^{(\ell)}_{\alpha\beta}=\delta_{k\ell}\
P^{(k)}_{\alpha\beta}. \]
\[ \bar P^{(k)}_{\alpha\beta}\ \bar P^{(\ell)}_{\alpha\beta}=\delta_{k\ell}\
\bar P^{(k)}_{\alpha\beta}. \]
{}From $\bar\qq_{\alpha\beta}\ff_{\alpha\beta}
= \ff_{\alpha\beta}\qq_{\alpha\beta}$ we have
\beq \bar P^{(k)}_{\alpha\beta}\ff_{\alpha\beta}
=\ff_{\alpha\beta}P^{(k)}_{\alpha\beta}\eeq

leading to

\beq \prod_{i=1}^{\mbox{\tiny rank \bf g}}
P^{\langle i,c_i\rangle}_{\alpha\beta}\bar
P^{(k)}_{\alpha\beta}\ff_{\alpha\beta}
=\ff_{\alpha\beta}\prod_{i=1}^{\mbox{\tiny rank \bf g}}
P^{\langle i,c_i\rangle}_{\alpha\beta}
P^{(k)}_{\alpha\beta} \eeq

This identity determines the orthogonal $\ff_{\alpha\beta}$ uniquely whenever
there are complete sets of one-dimensional projectors $P^a_{\alpha\beta}$ and
$\bar P^a_{\alpha\beta}$ among the
$\prod_{i=1}^{\mbox{\tiny rank \bf g }}\limits
P^{\langle i,c_i\rangle}_{\alpha\beta}
P^{(k)}_{\alpha\beta}$ and $\prod_{i=1}^{\mbox{\tiny rank \bf g}}\limits
P^{\langle i,c_i\rangle}_{\alpha\beta}\bar P^{(k)}_{\alpha\beta}$,
respectively. (The limit
$\lim_{h\to 0}\limits\ff=\eins^{\tn 2}$ fixes all signs uniquely.) This then
gives
us natural bases of the representation space of
$\varrho_\alpha\tn \varrho_\beta$ for
general $h$ and for $h\!=\!0$. Switching to the parameter $q=\;$e${}^h$ for
later
convenience, we denote the basis vectors by $|a;q\rangle$ and $|a;1\rangle$
with
orthogonality and completeness relations
\[ \langle a;q|a';q\rangle=\delta_{aa'} \]
\[ \sum_a |a;q\rangle\langle a;q|=\eins_\alpha\eins_\beta, \]
and
\[ P^a_{\alpha\beta}=|a;1\rangle\langle a;1|, \]
\[ \bar P^a_{\alpha\beta}=|a;q\rangle\langle a;q|. \]
Thus we finally obtain
\beq\label{F} \ff_{\alpha\beta}=\sum_a |a;q\rangle\langle a;1|. \eeq

{}From eq.(\ref{Phi}) it follows that
\beq\label{Phi2} \Phi=[\ff_{23}(\id\tn\Delta)(\ff)]^{-1}\ff_{12}(\Delta\tn\id)
(\ff) \eeq
or
\beq  \Phi_{\alpha\beta\gamma}
=[\ff_{\beta\gamma}(\id\tn\Delta)(\ff)_{\alpha\beta\gamma}]^{-1}
\ff_{\alpha\beta}(\Delta\tn\id)(\ff)_{\alpha\beta\gamma}. \eeq

\newpage
\section{Examples for F-matrices}
\subsection{su($n$)}
Let now $\g=$ su$(n)$, $\varrho_\alpha$ and $\varrho_\beta$ two copies of the
fundamental,
$n$-dimensional representation and $e_i$ $(i=1,\ldots,n)$ its orthonormal
standard basis
with dual $e^i$. In this case we find that we indeed obtain a complete set of
one-dimensional projectors by the method introduced above. For general (real,
positive)
$q$ we
find the following basis in terms of the product basis $e_i\tn e_j$:

\beq\label{subasis}\bgar{lrl}
\mbox{for } i=j:& |ii;q\rangle & =e_i\tn e_i \\
\mbox{for } i<j:& |ij_+;q\rangle &
=\frac{1}{\sqrt{q+q^{-1}}}\left(\sqrt{q}\,e_i\tn e_j
+\sqrt{q^{-1}}\,e_j\tn e_i\right)\\
\mbox{for } i<j:& |ij_-;q\rangle &
=\frac{1}{\sqrt{q+q^{-1}}}\left(\sqrt{q^{-1}}\,e_i\tn e_j
-\sqrt{q}\,e_j\tn e_i\right)\\
\ear\eeq

That implies

\beq\bgar{lrl}
\mbox{for } i=j:& \langle ii;1| & =e^i\tn e^i \\
\mbox{for } i<j:& \langle ij_+;1| &
=\frac{1}{\sqrt{2}}(e^i\tn e^j+e^j\tn e^i)\\
\mbox{for } i<j:& \langle ij_-;1| &
=\frac{1}{\sqrt{2}}(e^i\tn e^j-e^j\tn e^i)\\
\ear\eeq

and, finally,\footnote{We write $F$ for $\ff_{\alpha\beta}$ whenever we give
an explicit matrix
to emphasize that relations containing $\ff_{\alpha\beta}$ are valid for
general
$\varrho_\alpha$ and $\varrho_\beta$.}

\beq\bgar{ll} F= & \sum_i\limits (e_i\tn e_i)(e^i\tn e^i)\\
& + \sum_{i<j}\limits \frac{1}{\sqrt{2(q+q^{-1})}}
\left[\left(\sqrt{q}\,e_i\tn e_j+\sqrt{q^{-1}}\,e_j\tn e_i\right)
\left(e^i\tn e^j+e^j\tn e^i\right)\right.\\
& \left.+\left(\sqrt{q^{-1}}\,e_i\tn e_j-\sqrt{q}\,e_j\tn e_i\right)
\left(e^i\tn e^j-e^j\tn e^i\right)\right]. \\ \ear \eeq

For general representations of su($n$) we have one projector for every Young
diagram.
As far as the
eigenvalues of the Cartan generators uniquely fix a one-dimensional subspace
of an
irreducible representation it is guaranteed that our procedure will yield a
unique
F-matrix for an arbitrary pair of finite-dimensional representations
($\varrho_\alpha,\varrho_\beta$) of su($n$). If there are remaining
degeneracies one might
have to look for sufficiently natural additional input to fix a unique basis,
as will be
done below in the case of the defining representations of the Lie algebras of
the series
$B, C, D$.

\subsection{so($n$) and sp($n$)}

If we try to apply our method to a pair of fundamental representations of
so$(n)$ or
sp$(n)$ it turns out that there are still some degeneracies.
A closer look at the eigenspaces of the Cartan generators reveals beside
several
one- and two-dimensional subspaces (the latter
decompose in a symmetric and an antisymmetric one, as in the case of su$(n)$)
in this case also an
$n$-dimensional subspace. With the trace projector disposing us of one
dimension we are
left with $n\!-\!1$ dimensions to distribute among the symmetric and
antisymmetric
projectors.
Since this problem is increasing with the dimension of the Lie algebra's
defining
representation an inductive
procedure might help.

For so$(3)$ there is no degeneracy at all. In the case of so$(4)$ we are left
with a
two-dimensional antisymmetric subspace and in the case of sp$(4)$ with a
two-dimensional
symmetric one.

If one takes a look at the explicit expressions for the R-matrices of the
fundamental
representations of the three series $B, C, D$ [2] one notices their
characteristic
`onion-like'
structure. In the center of such an R-matrix one finds the R-matrices of all
the Lie
algebras of the same series with lower rank. So one can construct the Lie
algebras of
higher rank from those of lower rank in the same series by `adding coordinates
at the
outside'. This shows the way to a reasonably natural assumption that will fix
a unique
orthonormal
basis of the representation space of the tensor product of two defining
representations
of the Lie algebras of these series, and thus provide an F-matrix.

We assume that the symmetric and antisymmetric eigenvectors for the Lie
algebra of one
of these series of rank $r\!+\!1$ are the same as for the one of rank $r$ of
the same series
plus one additional symmetric and antisymmetric one, respectively, which are
determined by
the orthogonality and normalization requirement.

We denote by $n$
the dimension of the defining representation of a Lie algebra of one of the
series
$B, C, D$. For
$\varrho_\alpha$ and $\varrho_\beta$ we take two copies thereof. We
additionally use
$s=\frac n2$ (for $B$ and $D$) resp. $s=\frac n2+1$ (for $C$). As before we
denote the
standard basis vectors of the representation space by $e_i$ and their duals
by $e^i$. Again
we express our basis for the representation space of the tensor product of the
two
representations in terms of the product basis $e_i\tn e_j$.

In the case
$i\!+\!j\ne n\!+\!1$ we find the corresponding basis vectors of the same form
as
eq.(\ref{subasis}).
Using $\overline{\imath}=n\!+\!1\!-\!i$ and $\{k\}=q^k-q^{-k}$ we obtain by
the inductive
construction in the remaining cases the following basis vectors.\footnote{Thus
for
$q\!=\!1$ the expressions are to be replaced by their
respective limits.
Note that \mbox{$\lim_{q\!\to\! 1}\limits\frac{\{k\}}{\{m\}}=\frac km$} and
$\lim_{q\!\to\! 1}\limits \{k\}=0$.}
In terms of these the
F-matrices
are then given
by eq.(\ref{F}).

$D_s$:

\beq\bgar{rl}
\multicolumn{2}{l}{\mbox{for } k=1,\ldots,s\!-\!1:}\\
|n_{+(k)};q\rangle\;= & \sqrt{\frac{\{1\}^3}{\{2\}\{k\}\{k+1\}}}
\Bigg[\frac{\{k\}}{\{1\}}(q\,e_{s-k}\tn e_{\overline{s-k}}
+ q^{-1}\,e_{\overline{s-k}}\tn e_{s-k})\\&
-\sum_{i=0}^{k-1}\limits(q^{-i}\,e_{s-i}\tn e_{\overline{s-\imath}}
+ q^i\,e_{\overline{s-\imath}}\tn
e_{s-i})\Bigg]\\
\langle n_{+(k)};1|\;= & \sqrt{\frac{1}{2k(k+1)}}
\Bigg[k(e^{s-k}\tn e^{\overline{s-k}} + e^{\overline{s-k}}\tn e^{s-k})\\&
-\sum_{i=0}^{k-1}\limits(e^{s-i}\tn e^{\overline{s-\imath}}
+e^{\overline{s-\imath}}\tn
e^{s-i})\Bigg]\\
&\\
\multicolumn{2}{l}{\mbox{for } k=0,\ldots,s\!-\!1:}\\
|n_{-(k)};q\rangle\;= & \sqrt{\frac{\{1\}^3\{k-1\}\{k\}}{\{2\}\{2k-2\}\{2k\}}}
\Bigg[\frac{\{2k-2\}}{\{1\}\{k-1\}}(e_{s-k}\tn e_{\overline{s-k}}
- e_{\overline{s-k}}\tn e_{s-k})\\&
+\sum_{i=0}^{k-1}\limits(q^{-i}\,e_{s-i}\tn e_{\overline{s-\imath}}
+ q^i\,e_{\overline{s-\imath}}\tn
e_{s-i})\Bigg]\\
\langle n_{-(k)};1|\;= & \sqrt{\frac{1}{2}}(e^{s-k}\tn e^{\overline{s-k}}
-e^{\overline{s-k}}\tn e^{s-k})\\
&\\
\multicolumn{2}{l}{\mbox{and}}\\
|n_{\mbox{\tiny Tr}};q\rangle\;= & \sqrt{\frac{\{1\}\{s-1\}}{\{s\}\{2s-2\}}}
\sum_{i=0}^{s-1}\limits(q^{-i}\,e_{s-i}\tn e_{\overline{s-\imath}}
+ q^i\,e_{\overline{s-\imath}}\tn
e_{s-i})\\
\langle n_{\mbox{\tiny Tr}};1|\;= & \sqrt{\frac{1}{2s}}
\sum_{j=0}^{2s}\limits e^j\tn e^{\overline{\jmath}}\\
\ear\eeq

$B_{s-\frac 12}$:

\beq\bgar{rl}
\multicolumn{2}{l}{\mbox{for } k=\frac{1}{2},\frac{3}{2},\ldots,s\!-\!1:}\\
|n_{+(k)};q\rangle\;= & \sqrt{\frac{\{1\}^3}{\{2\}\{k\}\{k+1\}}}
\Bigg[\frac{\{k\}}{\{1\}}(q\,e_{s-k}\tn e_{\overline{s-k}}
+ q^{-1}\,e_{\overline{s-k}}\tn e_{s-k})\\&
-e_{s+\frac 12}\tn e_{s+\frac 12}
-\sum_{i=\frac{1}{2}}^{k-1}\limits(q^{-i}\,e_{s-i}\tn e_{\overline{s-\imath}}
+ q^i\,e_{\overline{s-\imath}}\tn
e_{s-i})\Bigg]\\
\langle n_{+(k)};1|\;= & \sqrt{\frac{1}{2k(k+1)}}
\Bigg[k(e^{s-k}\tn e^{\overline{s-k}} + e^{\overline{s-k}}\tn e^{s-k})\\&
-e^{s+\frac 12}\tn e^{s+\frac 12}
-\sum_{i=\frac{1}{2}}^{k-1}\limits(e^{s-i}\tn e^{\overline{s-\imath}}
+e^{\overline{s-\imath}}\tn
e^{s-i})\Bigg]\\
&\\
\multicolumn{2}{l}{\mbox{for } k=\frac{1}{2},\frac{3}{2},\ldots,s\!-\!1:}\\
|n_{-(k)};q\rangle\;= & \sqrt{\frac{\{1\}^3\{k-1\}\{k\}}{\{2\}\{2k-2\}\{2k\}}}
\Bigg[\frac{\{2k-2\}}{\{k-1\}\{1\}}(e_{s-k}\tn e_{\overline{s-k}}
- e_{\overline{s-k}}\tn e_{s-k})\\&+e_{s+\frac 12}\tn e_{s+\frac 12}
+\sum_{i=\frac{1}{2}}^{k-1}\limits(q^{-i}\,e_{s-i}\tn e_{\overline{s-\imath}}
+ q^i\,e_{\overline{s-\imath}}\tn
e_{s-i})\Bigg]\\
\langle n_{-(k)};1|\;= & \sqrt{\frac{1}2}(e^{s-k}\tn e^{\overline{s-k}}
-e^{\overline{s-k}}\tn e^{s-k})\\
&\\
\multicolumn{2}{l}{\mbox{and}}\\
|n_{\mbox{\tiny Tr}};q\rangle\;= & \sqrt{\frac{\{1\}\{s-1\}}{\{s\}\{2s-2\}}}
\Bigg[e_{s+\frac 12}\tn e_{s+\frac 12}\\&
+\sum_{i=\frac{1}{2}}^{s-1}\limits(q^{-i}\,e_{s-i}\tn e_{\overline{s-\imath}}
+ q^i\,e_{\overline{s-\imath}}\tn
e_{s-i})\Bigg]\\
\langle n_{\mbox{\tiny Tr}};1|\;= & \sqrt{\frac 1{2s}}
\sum_{j=0}^{2s}\limits e^j\tn e^{\overline{\jmath}}\\
\ear\eeq

$C_{s-1}$:

\beq\bgar{rl}
\multicolumn{2}{l}{\mbox{for } k=1,\ldots,s\!-\!1:}\\
|n_{+(k)};q\rangle\;= & \sqrt{\frac{\{1\}^3\{k\}\{k+1\}}{\{2\}\{2k\}\{2k+2\}}}
\Bigg[\frac{\{2k\}}{\{1\}\{k\}}(q\,e_{s-k}\tn e_{\overline{s-k}}
+ q^{-1}\,e_{\overline{s-k}}\tn e_{s-k})\\&
+\sum_{i=1}^{k-1}\limits(q^{-i}\,e_{s-i}\tn e_{\overline{s-\imath}}
- q^i\,e_{\overline{s-\imath}}\tn
e_{s-i})\Bigg]\\
\langle n_{+(k)};1|\;= & \sqrt{\frac{1}2}(e^{s-k}\tn e^{\overline{s-k}}
+e^{\overline{s-k}}\tn e^{s-k})\\
&\\
\multicolumn{2}{l}{\mbox{for } k=2,\ldots,s\!-\!1:}\\
|n_{-(k)};q\rangle\;= & \sqrt{\frac{\{1\}^3}{\{2\}\{k-1\}\{k\}}}
\Bigg[\frac{\{k-1\}}{\{1\}}(e_{s-k}\tn e_{\overline{s-k}}
- e_{\overline{s-k}}\tn e_{s-k})\\&
+\sum_{i=1}^{k-1}\limits(-q^{-i}\,e_{s-i}\tn e_{\overline{s-\imath}}
+ q^i\,e_{\overline{s-\imath}}\tn
e_{s-i})\Bigg]\\
\langle n_{-(k)};1|\;= & \sqrt{\frac{1}{2(k-1)k}}
\Bigg[(k\!-\!1)(e^{s-k}\tn e^{\overline{s-k}}
- e^{\overline{s-k}}\tn e^{s-k})\\&
+\sum_{i=1}^{k-1}\limits(-e^{s-i}\tn e^{\overline{s-\imath}}
+ e^{\overline{s-\imath}}\tn e^{s-i})\Bigg]\\
&\\
\multicolumn{2}{l}{\mbox{and}}\\
|n_{\mbox{\tiny Tr}};q\rangle\;= & \sqrt{\frac{\{1\}\{s\}}{\{s-1\}\{2s\}}}
\sum_{i=1}^{s-1}\limits(q^{-i}\,e_{s-i}\tn e_{\overline{s-\imath}}
- q^i\,e_{\overline{s-\imath}}\tn e_{s-i})\\
\langle n_{\mbox{\tiny Tr}};1|\;= & \sqrt{\frac{1}{2(s-1)}}
\sum_{i=1}^{s-1}\limits(e^{s-i}\tn e^{\overline{s-\imath}}
- e^{\overline{s-\imath}}\tn e^{s-i})\\
\ear\eeq

\newpage
\subsection{The crystal limit $q\!\to\! 0$}
The explicit expressions for the orthonormal basis vectors for general $q$,
$|a;q\rangle$,
now allow the
following observation:

For each basis vector $|a;q\rangle$ the limit $q\!\to\! 0$ exists and the set
of the
$|a;0\rangle$ is again an orthonormal basis of the representation space of the
respective
tensor product of representations.\footnote{Since in general
$\ff_{\alpha\beta}(q^{-1})=\ff_{\beta\alpha}(q)$ one can discuss the limit
$q\to\infty$ in
full analogy.}
In terms of the product basis $e_i\tn e_j$ it turns out to be particularly
simple:

\beq\bgar{rl}
|ii;0\rangle\;= & e_i\tn e_i \\
|ij_+;0\rangle\;= & e_j\tn e_i\\
|ij_-;0\rangle\;= & e_i\tn e_j\\
|n_{+(k)};0\rangle\;= &e_{\overline{s-k}}\tn e_{s-k}\\
|n_{-(k)};0\rangle\;= &\pm e_{s-k+1}\tn e_{\overline{s-k+1}}\qquad(+ \mbox{
for $B$ and $D$,
$-$ for $C$})\\
|n_{\mbox{\tiny Tr}};0\rangle\;= &e_1\tn e_n\\
\ear\eeq

In this limit the tensor product of two standard basis vectors is a basis
vector of an
irreducible component of the product representation. Thus transition from the
basis for
$q\!=\!1$ to the basis for $q\!=\!0$ means transition from a factorized basis
to a fully
reduced one, and the corresponding F-matrix' entries are therefore the
Clebsch-Gordan-coefficients.

This is in agreement with results obtained by Date, Jimbo, Miwa [4] and
Kashiwara [5].
Kashiwara calls the limit $q\!\to\! 0$ {\sl crystal limit} and the respective
bases {\sl crystal
bases}.

The existence of $\lim_{q\to 0}\limits\ff$ is not ensured by the
Drinfel'd-Kohno theorem
since in particular for the universal R-matices $\rr$ and $\bar\rr$ the limit
does not
exist. For unitary $\ff$ it is however evident from our construction that its
representations survive the limit which furthermore is smooth in every respect.

For two values $q,q'$ we will now write
\[ \ff^{[q'q]}=\ff'\ff^{-1} \]
such that
\[ \ff^{[q''q]}=\ff^{[q''q']}\ff^{[q'q]} \]
\[ \ff^{[q'q]}_{\alpha\beta}=\sum_a\limits |a;q'\rangle\langle a;q| \]
\[ \ff=\ff^{[q1]} \]
\[ \ff^{-1}=\ff^{[1q]}, \]
and, in particular,
\[ \ff^{[q1]}=\ff^{[q0]}\ff^{[01]}. \]

The representations of the coboundary of $\ff^{[01]}$ (with respect to the
``undeformed''
coproduct) can then be identified as the Racah coefficients.
\beq d\ff^{[01]}
=\ff^{[01]}_{23}(\id\tn\Delta)(\ff^{[01]})(\Delta\tn\id)(\ff^{[10]})
\ff^{[10]}_{12} \eeq
Note that this is what becomes of the (second) trivial coassociator
$\Phi_0\!=\!\eins^{\tn 3}$
of $U\g[[h]]$ under the twist with $\ff^{[q1]}$ (cf. eq. \ref{Phi}) in the
limit
$q\!\to\! 0$.
It is crucial that this
coassociator is {\sl not} compatible with the {\sl quasitriangular} structure
considered in
the Drinfel'd-Kohno theorem but with the ``classical'' {\sl triangular}
structure given by
the universal R-matrices $\rr_0\!=\!\eins^{\tn 2}$ and
$\bar\rr_0\!=\!\ff_{21}\ff^{-1}_{12}$. In
full analogy the expression for $\Phi$ given in eq.(\ref{Phi2}) is the
coboundary of
$\ff^{[q1]-1}\!=\!\ff^{[1q]}$ with respect to the ``deformed'' coproduct
\beq \bar d\ff^{[1q]}=\Phi.\eeq

It now turns out that the matrices $\ff^{[q0]}_{\alpha\beta}$ are of much
simpler form
than the matrices $\ff^{[q1]}_{\alpha\beta}$. E.g. for a pair of fundamental
representations
of su(2) it is
\beq F^{[q1]}=\left(
\bgar{cccc}
1&0&0&0\\
0&\frac{\sqrt {q^{}} + \sqrt{q^{-1}}}{\sqrt{2(q+q^{-1})}}
&\frac{\sqrt q - \sqrt{q^{-1}}}{\sqrt{2(q+q^{-1})}}&0\\
0&-\frac{\sqrt q - \sqrt{q^{-1}}}{\sqrt{2(q+q^{-1})}}
&\frac{\sqrt q + \sqrt{q^{-1}}}{\sqrt{2(q+q^{-1})}}&0\\
0&0&0&1\\
\ear\right)=\left(\bgar{cccc}1&0&0&0\\0&\frac{\sin\varphi+\cos\varphi}
{\sqrt 2}
&\frac{\sin\varphi-\cos\varphi}{\sqrt 2}&0\\
0&-\frac{\sin\varphi-\cos\varphi}{\sqrt 2}&\frac{\sin\varphi+\cos\varphi}
{\sqrt 2}&0\\0&0&0&1\\
\ear\right)
\eeq

\beq F^{[q0]}=\left(\bgar{cccc}
1&0&0&0\\0&\frac{\sqrt{q^{-1}}}{\sqrt{q+q^{-1}}}&
\frac{\sqrt q}{\sqrt{q+q^{-1}}}&0\\
0&-\frac{\sqrt q}{\sqrt{q+q^{-1}}}&\frac{\sqrt{q^{-1}}}{\sqrt{q+q^{-1}}}&0\\
0&0&0&1\\
\ear\right)=\left(\bgar{cccc}1&0&0&0\\0&\cos\varphi&\sin\varphi&0\\
0&-\sin\varphi&\cos\varphi&0\\0&0&0&1\\ \ear\right) \eeq

where $q=\tan\varphi$. In higher representations the fact that the form of
the
$\ff^{[q0]}_{\alpha\beta}$ is simpler will be more distinctive.

This is not so surprising if one takes into consideration that the quality of
$\ff^{[01]}_{\alpha\beta}$ to transform to a fully reduced basis simply means
that it
diagonalizes $\t_{\alpha\beta}$, and therefore $\qq_{\alpha\beta}$. From this
point of view
it makes good sense that $q\!=\!0$ is the natural reference point, rather than
$q\!=\!1$.
This is
because it lies in the very nature of tensor products that there is no inherent
information
that would allow an identification as a ``composed'' object -- as opposed to an
``elementary'' one. Thus from the perspective of the product representation the
factorized
basis ($q\!=\!1$) is not distinguished in any way from others. The only one for
which this
is the case is the
fully reduced one ($q\!=\!0$).

\section{Conclusion}
Suppose you have two matrix representations of a compact simple Lie group with
orthonormal bases
chosen such
that the Cartan generators are diagonal.
There are two natural choices for orthonormal
bases for
their tensor product, the product basis and the fully reduced basis with
respect to which
the Cartan generators are still diagonal. Transition between them can be
thought of as
rotations in the simultaneous eigenspaces of the Cartan generators. This way
one can
define a one-parameter family of bases, in all of which the Cartan generators
are diagonal,
assigning the parameter value $0$ to the fully reduced basis and the value $1$
to the
factorized basis such that for every parameter value the basis is orthonormal
and the
dependence on the parameter is smooth. Then this parameter can be chosen as
Drinfel'd's
``deformation'' parameter $q\!=\!\;$e${}^h$ and the one-parameter family of
bases is
given by $|a;q\rangle$.

The Drinfel'd-Kohno theorem shows that quantum groups are equivalent to
compact simple
Lie groups with an additional structure. It turns out that the additional
structure amounts
to the information about the degree of factorization of the natural basis of
the
tensorproduct of two matrix representations, expressed in terms of the
``deformation''
parameter $q$. In view of this it might be more appropriately called a
factorization
parameter.

Results of this work are applied in [6]. A further application by the author
is in
progress [7].

{\sl I would like to thank Julius Wess for encouragement and for the
opportunity
to do this work under
his supervision, John C. Baez, John F. Cornwell, Gaetano Fiore, Dale
Husemoller, Shahn Majid
and Peter Schupp for valuable
suggestions and Kristin F\"orger for helping me understand a crucial part of
Drinfel'd's
theory.}

\end{document}